\documentstyle[12pt,leqno]{article}

\font\twlgot =eufm10 scaled \magstep1
\font\egtgot =eufm8
\font\sevgot =eufm7

\font\twlmsb =msbm10 scaled \magstep1
\font\egtmsb =msbm8
\font\sevmsb =msbm7

\newfam\gotfam
\def\pgot{\fam\gotfam\twlgot}
\textfont\gotfam\twlgot
\scriptfont\gotfam\egtgot
\scriptscriptfont\gotfam\sevgot
\def\got{\protect\pgot}

\newfam\msbfam
\textfont\msbfam\twlmsb
\scriptfont\msbfam\egtmsb
\scriptscriptfont\msbfam\sevmsb
\def\Bbb{\protect\pBbb}
\def\pBbb{\relax\ifmmode\expandafter\Bb\else\typeout{You cann't use
Bbb in text mode}\fi}
\def\Bb #1{{\fam\msbfam\relax#1}}

\newcommand{\gO}{{\got O}}
\newcommand{\gQ}{{\got Q}}
\newcommand{\gU}{{\got U}}
\newcommand{\gE}{{\got E}}

\def\thebibliography#1{\centerline{ \sc References}\list
  {[\arabic{enumi}]}{\settowidth\labelwidth{#1}\leftmargin\labelwidth
    \advance\leftmargin\labelsep
    \usecounter{enumi}}
    \def\newblock{\hskip .11em plus .33em minus .07em}
    \sloppy\clubpenalty4000\widowpenalty4000
    \sfcode`\.=1000\relax}

\def\op#1{\mathop{\fam0 #1}\limits}

\newcommand{\id}{{\rm Id\,}}

\newcommand{\di}{{\rm dim\,}}

\newcommand{\nm}[1]{\mid {#1}\mid}

\newcommand{\beq}{\begin{equation}}
\newcommand{\eeq}{\end{equation}}
\newcommand{\ben}{\begin{eqnarray}}
\newcommand{\een}{\end{eqnarray}}
\newcommand{\be}{\begin{eqnarray*}}
\newcommand{\ee}{\end{eqnarray*}}
\newcommand{\bea}{\begin{eqalph}}
\newcommand{\eea}{\end{eqalph}}

\newcommand{\cO}{{\cal O}}

\newcommand{\cQ}{{\cal Q}}
\newcommand{\cE}{{\cal E}}

\newcommand{\al}{\alpha}

\newcommand{\dl}{\delta}
\newcommand{\la}{\lambda}
\newcommand{\La}{\Lambda}
\newcommand{\f}{\phi}

\newcommand{\m}{\mu}

\newcommand{\G}{\Gamma}
\newcommand{\th}{\theta}

\newcommand{\si}{\sigma}

\newcommand{\w}{\wedge}

\newcommand{\ol}{\overline}
\newcommand{\dr}{\partial}
\newcommand{\ar}{\op\longrightarrow}
\newcommand{\ot}{\otimes}

\newcounter{eqalph}
\newcounter{equationa}
\newcounter{example}
\newcounter{remark}
\newcounter{theorem}
\newcounter{proposition}
\newcounter{lemma}
\newcounter{corollary}
\newcounter{definition}

\def\thedefinition{\arabic{definition}}

\newenvironment{proof}{\noindent{\it Proof.}}{\hfill $\Box$
\medskip }
\newenvironment{rem}{\refstepcounter{definition} \medskip \noindent{\bf Remark
\thedefinition}.}{ \medskip }
\newenvironment{theo}{\refstepcounter{definition} \medskip\noindent{\bf
Theorem \thedefinition.}\sl}{\medskip }
\newenvironment{prop}{\refstepcounter{definition} \medskip\noindent{\bf
Proposition \thedefinition.}\sl}{\medskip }
\newenvironment{lem}{\refstepcounter{definition} \medskip\noindent{\bf  Lemma
\thedefinition.}\sl }{\medskip }
\newenvironment{cor}{\refstepcounter{definition} \medskip\noindent{\bf 
Corollary \thedefinition.}\sl }{\medskip }

\newenvironment{eqalph}{\stepcounter{equation}
\setcounter{equationa}{\value{equation}}
\setcounter{equation}{0}

\begin{eqnarray}}{\end{eqnarray}\setcounter{equation}{\value{equationa}}}

\newcommand{\subs}[1]{\bigskip\centerline{\sc #1}\bigskip}

\begin{document}
\hbox{}

\begin{center}

{ \Large \bf Cohomology of the variational bicomplex on
the infinite order jet space}
\bigskip\bigskip

{\sc GIOVANNI GIACHETTA, LUIGI MANGIAROTTI, AND GENNADI SARDANASHVILY}
\bigskip
\bigskip

\end{center}

{\small
\noindent
{\sc Abstract.} 
We obtain the cohomology of the variational bicomplex on the infinite
order jet space of a smooth fiber bundle in the class of exterior forms of
finite jet order. This provides a solution of the global
inverse problem of the calculus of variations of finite order on fiber bundles.
\bigskip

\noindent
2000 {\it Mathematics subject classification.} Primary 58A20, 58E30; Secondary
55N30.
\medskip

\noindent
{\it Key words and phrases}.  Jet manifolds, infinite order
jet space, differential graded algebra, variational
complex, cohomology of sheaves.
}
\bigskip

\subs{1. Introduction}

Let $Y\to X$ be a smooth fiber bundle. We obtain cohomology of the 
variational bicomplex on the infinite
order jet space 
$J^\infty Y$ of $Y\to X$ in the class of exterior forms of finite jet order.
This is cohomology of the vertical differential
$d_V$, the horizontal (or total) differential $d_H$ and the variational
operator $\dl$.

The two differential calculus  of exterior
forms $\cO^*_\infty$ and $\cQ^*_\infty$ are usually considered on
$J^\infty Y$.  The $\cO^*_\infty$ is the direct
limit of graded differential algebras of exterior forms on finite
order jet manifolds. 
Its cohomology,  except de Rham cohomology and a particular
result of
\cite{vin} on $\dl$-cohomology, remains unknown. At the same time,
$\cO^*_\infty$ is most interesting for  applications because it consists of
exterior forms on finite order jet manifolds. The
$\cQ^*_\infty$ is the structure algebra  of the sheaf of germs of exterior
forms on finite order jet manifolds. There is the $\Bbb R$-algebra monomorphism
$\cO^*_\infty \to \cQ^*_\infty$. The $d_H$- and $\dl$-cohomology of
$\cQ^*_\infty$ has been investigated in \cite{ander80,tak2}. 
Due to Lemma \ref{lmp03} below, 
we simplify this investigation
and complete it by the study of $d_V$-cohomology of $\cQ^*_\infty$. We
prove that the graded differential algebra 
$\cO^*_\infty$ has the same $d_H$- and $\dl$-cohomology as $\cQ^*_\infty$
(see Theorem \ref{am11} below). This provides a solution of the global
inverse problem of the calculus of variations in the class of exterior forms
of finite jet order. 

Note that the local exactness of the calculus of variations
has been proved in the class of exterior forms of finite order by use of
homotopy operators which do not minimize the order of Lagrangians (see, e.g.,
\cite{olver,tul}). The infinite variational complex of such
exterior forms on $J^\infty Y$ has been studied by many authors (see,
e.g., \cite{bau,book,olver,tak1,tul}). However, these forms on $J^\infty Y$
fail to constitute a sheaf. Therefore, the cohomology obstruction to the
exactness of the calculus of variations has been obtained in the class of
exterior forms of locally finite order which make up the above
mentioned algebra $\cQ^*_\infty$
\cite{ander,tak2}. 
A solution of the global inverse problem in the calculus of variations in
the class of exterior forms of a fixed jet order has been suggested in
\cite{ander80} by a computation of cohomology of the fixed
order variational sequence (see
\cite{kru90,vit} for another variant of such a variational sequence).
The key point of this computation lies in the local exactness
of the finite order variational sequence which however requires
rather sophisticated {\it ad hoc} technique in order to be reproduced (see also
\cite{kru98}). We show that  the obstruction to
the exactness of the finite order calculus of variations is the same as 
for exterior forms of locally finite order, without minimizing an order
of Lagrangians.
The main point for applications is that this obstruction is given by closed
forms on the fiber bundle
$Y$, and is of first order. 

The article is organized as follows. In Section 2, the differential
calculus $\cO^*_\infty$ and $\cQ^*_\infty$ on $J^\infty Y$ are introduced  in
an algebraic way. In Section 3, the variational bicomplex on $J^\infty Y$ is
set. Section 4 is devoted to cohomology
of the differential calculus
$\cQ^*_\infty$ on
$J^\infty Y$. In Section 5, the isomorphism of $d_H$- and $\dl$-cohomology of
$\cQ^*_\infty$ to that of
$\cO^*_\infty$ is proved. In Sections 6, a solution of the global inverse
problem in the calculus of variations in different classes of exterior forms is
provided. 

\subs{2. The differential calculus on $J^\infty Y$}

Smooth manifolds throughout are assumed to be
real, finite-dimensional, Hausdorff,
paracompact, and connected. Put further dim$X=n\geq 1$. We follow the
standard terminology of jet formalism \cite{book,man,sard}.

Recall that the infinite order jet space $J^\infty Y$ of a smooth fiber bundle
$Y\to X$ is defined as a projective limit
$(J^\infty Y,\pi^\infty_r)$ of the inverse system
\beq
X\op\longleftarrow^\pi Y\op\longleftarrow^{\pi^1_0}\cdots \longleftarrow
J^{r-1}Y \op\longleftarrow^{\pi^r_{r-1}} J^rY\longleftarrow\cdots \label{t1}
\eeq
of finite order jet manifolds $J^rY$ of $Y\to X$, where $\pi^r_{r-1}$ are
affine bundles. Bearing in mind Borel's theorem, one can say that 
$J^\infty Y$ consists of the equivalence classes of sections of $Y\to X$
identified by their Taylor series at points of $X$.
Endowed with the projective limit topology,
$J^\infty Y$ is a paracompact Fr\'echet manifold \cite{tak2}.
A bundle coordinate atlas
$\{U_Y,(x^\la,y^i)\}$ of $Y\to X$ yields the manifold
coordinate atlas
\be
\{(\pi^\infty_0)^{-1}(U_Y), (x^\la, y^i_\La)\}, \qquad 0\leq|\La|,
\ee
 of $J^\infty
Y$, together with the transition functions  
\beq
{y'}^i_{\la+\La}=\frac{\dr x^\m}{\dr x'^\la}d_\m y'^i_\La, \label{55.21}
\eeq
where $\La=(\la_k\ldots\la_1)$, $\la+\La=(\la\la_k\ldots\la_1)$ are
multi-indices and
$d_\la$ denotes the total derivative 
\be
d_\la = \dr_\la + \op\sum_{|\La|\geq 0} y^i_{\la+\La}\dr_i^\La.
\ee

With the inverse system (\ref{t1}), one has
the direct system 
\be
\cO^*(X)\op\longrightarrow^{\pi^*} \cO^*_0 
\op\longrightarrow^{\pi^1_0{}^*} \cO_1^*
\op\longrightarrow^{\pi^2_1{}^*} \cdots \op\longrightarrow^{\pi^r_{r-1}{}^*}
 \cO_r^* \longrightarrow\cdots 
\ee
of graded differential $\Bbb R$-algebras $\cO^*_r$ of exterior forms on finite
order jet manifolds $J^rY$, where $\pi^r_{r-1}{}^*$ are pull-back
monomorphisms. The direct limit
 of this direct system is the above mentioned graded
differential $\Bbb R$-algebra $(\cO^*_\infty,\pi^{\infty *}_r)$ 
of exterior forms on  finite order jet manifolds modulo the pull-back
identification. The $\cO^*_\infty$ is a
differential calculus over the $\Bbb R$-ring $\cO^0_\infty$ of continuous real
functions on
$J^\infty Y$ which are the pull-back of smooth real functions on
finite order jet manifolds by surjections $\pi^\infty_r$. 
Passing to the direct limit of the de Rham complexes of exterior forms
on finite order
jet manifolds, de Rham cohomology of the graded
differential algebra
$\cO^*_\infty$ has been found, and coincides with de Rham cohomology
of the fiber bundle $Y$
\cite{ander,bau}. However, this is not a way of studying other cohomology of
the algebra $\cO^*_\infty$.  

To solve this problem, let us enlarge $\cO^0_\infty$ to the $\Bbb
R$-ring
$\cQ^0_\infty$ of continuous real functions on $J^\infty Y$ such that, given
$f\in
\cQ^0_\infty$ and any point $q\in J^\infty Y$, there exists a neighborhood
of $q$ where $f$ coincides with the pull-back of a smooth function on some
finite order jet manifold.  The reason lies in 
the fact that the paracompact space
$J^\infty Y$ admits a partition of unity by elements of the ring
$\cQ^0_\infty$ \cite{tak2}. Therefore, sheaves of
$\cQ^0_\infty$-modules on
$J^\infty Y$ are fine and, consequently, acyclic. Then, the 
abstract de Rham theorem on cohomology of a sheaf resolution can be called
into play. 

\begin{rem}
Throughout, we follow the terminology of
\cite{hir} where by a sheaf $S$ over a topological space $Z$ is meant a sheaf
bundle $S\to Z$. Accordingly, $\G(S)$ denotes the canonical presheaf of
sections of the sheaf $S$, and 
$\G(Z,S)$ is the group of global sections of $S$. All sheaves below are
ringed spaces, but we omit this terminology if there is no danger of
confusion.
\end{rem}

Let us define a differential calculus over the ring $\cQ^0_\infty$.
Let $\gO^*_r$ be a sheaf
of germs of exterior forms on the $r$-order jet manifold $J^rY$ and 
$\G(\gO^*_r)$ its canonical presheaf.  There is the direct system of canonical
presheaves
\be
\G(\gO^*_X)\op\longrightarrow^{\pi^*} \G(\gO^*_0) 
\op\longrightarrow^{\pi^1_0{}^*} \G(\gO_1^*)
\op\longrightarrow^{\pi^2_1{}^*} \cdots \op\longrightarrow^{\pi^r_{r-1}{}^*}
 \G(\gO_r^*) \longrightarrow\cdots, 
\ee
where $\pi^r_{r-1}{}^*$ are pull-back monomorphisms
with respect to open surjections 
$\pi^r_{r-1}$. Its direct
limit $\gO^*_\infty$ 
is a presheaf of graded differential
$\Bbb R$-algebras on
$J^\infty Y$. Let $\gQ^*_\infty$ be a sheaf constructed from 
$\gO^*_\infty$ and $\G(\gQ^*_\infty)$ its canonical presheaf. There is the 
$\Bbb R$-algebra monomorphism of presheaves 
 $\gO^*_\infty
\to\G(\gQ^*_\infty)$. 
The structure algebra 
$\cQ^*_\infty=\G(J^\infty Y,\gQ^*_\infty)$ of
the sheaf $\gQ^*_\infty$  
is a desired differential calculus over the $\Bbb
R$-ring $\cQ^*_\infty$. 

For short, we
agree to call elements of $\cQ^*_\infty$ the
exterior forms on
$J^\infty Y$.  Restricted to a
coordinate chart
$(\pi^\infty_0)^{-1}(U_Y)$ of $J^\infty Y$, they
can be written in a coordinate form, where horizontal forms
$\{dx^\la\}$ and contact 1-forms
$\{\th^i_\La=dy^i_\La -y^i_{\la+\La}dx^\la\}$ constitute the set of
generators of the algebra
$\cQ^*_\infty$. 
There is the canonical splitting
\be
\cQ^*_\infty =\op\oplus_{k,s}\cQ^{k,s}_\infty, \qquad 0\leq k, \qquad
0\leq s\leq n,
\ee
of $\cQ^*_\infty$ into $\cQ^0_\infty$-modules $\cQ^{k,s}_\infty$
of $k$-contact and $s$-horizontal forms, together with the corresponding
projections
\be
h_k:\cQ^*_\infty\to \cQ^{k,*}_\infty, \quad 0\leq k, \qquad
h^s:\cQ^*_\infty\to \cQ^{*,s}_\infty, \quad 0\leq s
\leq n.
\ee 
Accordingly, the
exterior differential on $\cQ_\infty^*$ is
decomposed into the sum $d=d_H+d_V$ of horizontal and vertical
differentials such that
\be
&& d_H\circ h_k=h_k\circ d\circ h_k, \qquad d_H(\f)=
dx^\la\w d_\la(\f), \qquad \f\in\cQ^*_\infty,\\ 
&& d_V \circ h^s=h^s\circ d\circ h^s, \qquad
d_V(\f)=\th^i_\La \w \dr_\La^i\f.
\ee

\subs{3. The variational bicomplex}

Being nilpotent, the
differentials $d_V$ and $d_H$ provide the natural bicomplex
$\{\gQ^{k,m}_\infty\}$ of  the sheaf
$\gQ^*_\infty$ on $J^\infty Y$. To complete it to the
variational bicomplex, one defines the projection $\Bbb R$-module
endomorphism 
\be
&& \tau=\op\sum_{k>0} \frac1k\ol\tau\circ h_k\circ h^n, \\ 
&&\ol\tau(\f)
=(-1)^{\nm\La}\th^i\w [d_\La(\dr^\La_i\rfloor\f)], \qquad 0\leq\nm\La,
\qquad \f\in \G(\gQ^{>0,n}_\infty),
\ee
of $\gQ^*_\infty$ such that
\be
\tau\circ d_H=0, \qquad  \tau\circ d\circ \tau -\tau\circ d=0.
\ee
Introduced on elements of the presheaf $\gO^*_\infty$ 
(see, e.g., \cite{bau,book,tul}), this endomorphism is induced on the
sheaf $\gQ^*_\infty$ and its structure algebra
$\cQ^*_\infty$. Put
\be
\gE_k=\tau(\gQ^{k,n}_\infty), \qquad E_k=\tau(\cQ^{k,n}_\infty), \qquad k>0.
\ee
Since
$\tau$ is a projection operator, we have isomorphisms 
\be
\G(\gE_k)=\tau(\G(\gQ^{k,n}_\infty)), \qquad E_k=\G(J^\infty Y,\gE_k).
\ee
The variational operator on $\gQ^{*,n}_\infty$ is defined as the
morphism $\dl=\tau\circ d$. 
It is nilpotent, and obeys the relation 
\beq
\dl\circ\tau-\tau\circ d=0. \label{am13}
\eeq

Let $\Bbb R$ and  $\gO^*_X$ denote the constant sheaf
on
$J^\infty Y$ and the sheaf of exterior forms on $X$, respectively. The
operators $d_V$,
$d_H$,
$\tau$ and $\dl$ give the following variational bicomplex of
sheaves of differential forms on $J^\infty Y$:
\beq
\begin{array}{ccccrlcrlccrlccrlcrl}
& & & & _{d_V} & \put(0,-7){\vector(0,1){14}} & & _{d_V} &
\put(0,-7){\vector(0,1){14}} & &  & _{d_V} &
\put(0,-7){\vector(0,1){14}} & & &  _{d_V} &
\put(0,-7){\vector(0,1){14}}& & _{-\dl} & \put(0,-7){\vector(0,1){14}} \\ 
 &  & 0 & \to & &\gQ^{k,0}_\infty &\ar^{d_H} & & \gQ^{k,1}_\infty &
\ar^{d_H} &\cdots  & & \gQ^{k,m}_\infty &\ar^{d_H} &\cdots & &
\gQ^{k,n}_\infty &\ar^\tau &  & \gE_k\to  0\\  
 & &  &  & & \vdots & & & \vdots  & & & 
&\vdots  & & & &
\vdots & &   & \vdots \\ 
& & & & _{d_V} & \put(0,-7){\vector(0,1){14}} & & _{d_V} &
\put(0,-7){\vector(0,1){14}} & &  & _{d_V} &
\put(0,-7){\vector(0,1){14}} & & &  _{d_V} &
\put(0,-7){\vector(0,1){14}}& & _{-\dl} & \put(0,-7){\vector(0,1){14}} \\ 
 &  & 0 & \to & &\gQ^{1,0}_\infty &\ar^{d_H} & & \gQ^{1,1}_\infty &
\ar^{d_H} &\cdots  & & \gQ^{1,m}_\infty &\ar^{d_H} &\cdots & &
\gQ^{1,n}_\infty &\ar^\tau &  & \gE_1\to  0\\  
& & & & _{d_V} &\put(0,-7){\vector(0,1){14}} & & _{d_V} &
\put(0,-7){\vector(0,1){14}} & & &  _{d_V}
 & \put(0,-7){\vector(0,1){14}} & &  & _{d_V} & \put(0,-7){\vector(0,1){14}}
 & & _{-\dl} & \put(0,-7){\vector(0,1){14}} \\
0 & \to & \Bbb R & \to & & \gQ^0_\infty &\ar^{d_H} & & \gQ^{0,1}_\infty &
\ar^{d_H} &\cdots  & &
\gQ^{0,m}_\infty & \ar^{d_H} & \cdots & &
\gQ^{0,n}_\infty & \equiv &  & \gQ^{0,n}_\infty \\
& & & & _{\pi^{\infty*}}& \put(0,-7){\vector(0,1){14}} & & _{\pi^{\infty*}} &
\put(0,-7){\vector(0,1){14}} & & &  _{\pi^{\infty*}}
 & \put(0,-7){\vector(0,1){14}} & &  & _{\pi^{\infty*}} &
\put(0,-7){\vector(0,1){14}} & &  & \\
0 & \to & \Bbb R & \to & & \gO^0_X &\ar^d & & \gO^1_X &
\ar^d &\cdots  & &
\gO^m_X & \ar^d & \cdots & &
\gO^n_X & \ar^d & 0 &  \\
& & & & &\put(0,-5){\vector(0,1){10}} & & &
\put(0,-5){\vector(0,1){10}} & & & 
 & \put(0,-5){\vector(0,1){10}} & & &   &
\put(0,-5){\vector(0,1){10}} & &  & \\
& & & & &0 & &  & 0 & & & & 0 & & & & 0 & &  & 
\end{array}
\label{7}
\eeq
The second row and the last column of this bicomplex form the 
variational complex
\beq
0\to\Bbb R\to \gQ^0_\infty \ar^{d_H}\gQ^{0,1}_\infty\ar^{d_H}\cdots  
\op\longrightarrow^{d_H} 
\gQ^{0,n}_\infty  \op\longrightarrow^\dl \gE_1 
\op\longrightarrow^\dl 
\gE_2 \longrightarrow \cdots\, . \label{tams1}
\eeq
The corresponding variational bicomplexes $\{\cQ^*_\infty,E_k\}$ and
$\{\cO^*_\infty, \ol E_k\}$ of the differential calculus $\cQ^*_\infty$ and 
$\cO^*_\infty$ take place.

There are the well-known statements summarized usually as
the algebraic Poincar\'e lemma (see, e.g., \cite{olver,tul}). 

\begin{lem} \label{am12} 
If $Y$ is a contractible fiber bundle $\Bbb R^{n+p}\to\Bbb R^n$, the
variational bicomplex $\{\cO^*_\infty, \ol E_k\}$ of the graded differential
algebra $\cO^*_\infty$ is exact.
\end{lem}

It follows that the variational
bicomplex of sheaves (\ref{7}) is exact for any smooth fiber bundle $Y\to X$.
Moreover, all sheaves
$\gQ^{k,m}$ in this bicomplex are fine, and so are the sheaves $\gE_k$ in
accordance with the following lemma.

\begin{lem} \label{lmp03} 
Sheaves $\gE_k$, $k>0$, are fine.
\end{lem}

\begin{proof}
Though $\Bbb R$-modules $E_{k>1}$ fail to be
$\cQ^0_\infty$-modules \cite{tul}, one can use the fact that the sheaves
$\gE_{k>0}$ are projections $\tau(\gQ^{k,n}_\infty)$ of sheaves of
$\cQ^0_\infty$-modules. Let $\gU =\{U_i\}_{i\in I}$
be a 
locally finite open covering  of
$J^\infty Y$ and $\{f_i\in\cQ^0_\infty\}$ the associated partition of unity. 
For any open subset $U\subset J^\infty Y$ and any section
$\varphi$ of 
the sheaf $\gQ^{k,n}_\infty$ over $U$, let us put
$h_i(\varphi)=f_i\varphi$. Then, $\{h_i\}$ provide a family of endomorphisms
of the sheaf $\gQ^{k,n}_\infty$, required for $\gQ^{k,n}_\infty$ to be fine.
Endomorphisms $h_i$ of $\gQ^{k,n}_\infty$ also yield the $\Bbb R$-module
endomorphisms 
\be
\ol h_i=\tau\circ h_i: \gE_k\ar^{\rm in} \gQ^{k,n}_\infty \ar^{h_i}
\gQ^{k,n}_\infty \ar^\tau \gE_k
\ee
of the sheaves $\gE_k$.
They possess the properties
required for $\gE_k$ to be a fine sheaf. Indeed, for each $i\in I$, there is a
closed set ${\rm supp}\,f_i\subset U_i$ such that $\ol h_i$ is zero outside
this set, while the sum $\op\sum_{i\in I}\ol h_i$ is the identity morphism.
\end{proof}

This Lemma simplify essentially our cohomology computation of the
variational bicomplex in comparison with that in \cite{ander80,tak2}.  
Since all sheaves except $\Bbb R$ and $\pi^{\infty*}\gO^*_X$ in the bicomplex
(\ref{7}) are fine, the abstract de Rham theorem 
(\cite{hir}, Theorem 2.12.1) can be applied to columns
and rows of this bicomplex in a straightforward way. We will quote the
following variant of this theorem.

\begin{theo} \label{+132} 
Let 
\beq
0\to S\ar^h S_0\ar^{h^0} S_1\ar^{h^1}\cdots\ar^{h^{p-1}} S_p\ar^{h^p}
S_{p+1}, \qquad p>1, \label{+131'}
\eeq
be an exact sequence of sheaves on a paracompact topological space $Z$, where
the sheaves $S_p$ and $S_{p+1}$ are not necessarily acyclic, and let 
\beq
0\to \G(Z,S)\ar^{h_*} \G(Z,S_0)\ar^{h^0_*}
\G(Z,S_1)\ar^{h^1_*}\cdots\ar^{h^{p-1}_*} \G(Z,S_p)\ar^{h^p_*}
\G(Z,S_{p+1}) \label{+130}
\eeq
be the corresponding cochain complex
of structure groups of these sheaves.
The $q$-cohomology groups of the
cochain complex (\ref{+130}) for $0\leq q\leq p$ are
isomorphic to the cohomology groups $H^q(Z,S)$ of $Z$ with coefficients in the
sheaf $S$. 
\end{theo}

The $d_H$- and $\dl$-cohomology of the differential calculus $\cQ^*_\infty$ 
on $J^\infty Y$ has been found in \cite{tak2}. Its
$\dl$-cohomology at terms $\cQ^{0,n}_\infty$ and
$E_1$ has been also obtained in \cite{ander80} (several statements without
proof were announced in \cite{ander}). We recover this cohomology
in a short way due to Lemma \ref{lmp03}, and complete it by the
$d_V$-cohomology of
$\cQ^*_\infty$ corresponding to columns of the bicomplex (\ref{7}).

\subs{4. Cohomology of $\cQ^*_\infty$}

We start from the following facts.

\begin{lem} \label{20jpa} 
There is an
isomorphism 
\beq
H^*(J^\infty Y,\Bbb R)= H^*(Y,\Bbb R)=H^*(Y) \label{lmp80}
\eeq
between cohomology $H^*(J^\infty Y,\Bbb R)$ of $J^\infty Y$ with
coefficients in the constant sheaf $\Bbb R$, that $H^*(Y,\Bbb R)$ of $Y$, and
de Rham cohomology $H^*(Y)$ of $Y$. 
\end{lem}

\begin{proof}
A fiber bundle $Y$ is a strong deformation retract of $J^\infty Y$. Then,
the first isomorphism in (\ref{lmp80}) follows from the Vietoris--Begle
theorem (\cite{bred}, Theorem 11.4; \cite{kash}, Corollary 2.7.7), while the
second one is a consequence of  the well-known de Rham theorem.
\end{proof}

\begin{lem} \label{20jpa'} 
There is an isomorphism
\beq
H^*(J^\infty Y,\pi^{\infty *}\gO^m_X)=H^*(Y,\pi^*\gO^m_X) \label{7jp'}
\eeq
between cohomology $H^*(J^\infty Y,\pi^{\infty *}\gO^m_X)$ of $J^\infty Y$
with coefficients in the pull-back sheaf $\pi^{\infty *}\gO^m_X$ and that of
$Y$ with coefficients in the sheaf $\pi^*\gO^m_X$.
\end{lem}

\begin{proof}
The isomorphism (\ref{7jp'}) also follows from the facts that
$Y$ is a strong
deformation retract of $J^\infty Y$ and that $\pi^{\infty*}\gO^m_X$ is the
pull-back onto $J^\infty Y$ of the sheaf $\pi^*\gO^m_X$ on $Y$ 
(\cite{kash}, Corollary 2.7.7). 
\end{proof}

\begin{rem}
Lemma \ref{20jpa} and Lemma \ref{20jpa'} are corollary of Lemma \ref{am12} as
follows. Let us consider the open surjection $\pi^\infty_0:J^\infty Y\to Y$
and the direct images $\pi^\infty_{0*}\Bbb R$ and
$\pi^\infty_{0*}(\pi^{\infty *}\gO^m_X)$ of sheaves $\Bbb R$ and $\pi^{\infty
*}\gO^m_X$ on $J^\infty Y$. They are isomorphic to the sheaves
$\Bbb R$ and $\pi^*\gO^m_X$ on $Y$, respectively. Lemma \ref{am12} shows that,
every point
$y\in Y$ has a base of open neighbourhoods $\{U_y\}$ whose inverse images 
$(\pi^\infty_0)^{-1}(U_y)$
are acyclic for the sheaves $\Bbb R$ and $\pi^{\infty
*}\gO^m_X$. Then, a weak version of the Leray theorem states the cohomology
isomorphisms (\ref{lmp80}) and (\ref{7jp'}) \cite{god}. Moreover, since other
sheaves in the bicomplex (\ref{7}) are acyclic on
$(\pi^\infty_0)^{-1}(U_y)$ and the bicomplexes of their sections over
$(\pi^\infty_0)^{-1}(U_y)$ are exact, we have the exact direct image
$\{\pi^\infty_{0*}\gQ^*_\infty,\pi^\infty_{0*}\gE_k\}$ on $Y$ of the bicomplex
(\ref{7}), whose rows and columns are
resolutions of sheaves on
$Y$.  Due to the
$\Bbb R$-algebra isomorphism
$\cQ^*_\infty=\G(Y,\pi^\infty_{0*}\gQ^*_\infty)$, one can study cohomology of
the graded differential $\Bbb R$-algebra
$\cQ^*_\infty$ by use of this variational bicomplex on $Y$. In particular, it
follows that cohomology of $\cQ^*_\infty$ of degree $q>\di Y$ vanishes.
\end{rem}

Turn to de Rham cohomology of the algebra $\cQ^*_\infty$. 
Let us consider the de Rham complex of sheaves 
 \beq
0\to \Bbb R\to
\gQ^0_\infty\op\longrightarrow^d\gQ^1_\infty\op\longrightarrow^d
\cdots
\label{lmp71}
 \eeq
on $J^\infty Y$ and the de Rham complex of their structure algebras
 \beq
0\to \Bbb R\to
\cQ^0_\infty\op\longrightarrow^d\cQ^1_\infty\op\longrightarrow^d
\cdots\, .
\label{5.13'}
\eeq

\begin{prop} \label{38jp}  There is an isomorphism
\be
H^*(\cQ^*_\infty)=H^*(Y) 
\ee
of de Rham cohomology $H^*(\cQ^*_\infty)$ 
of the graded differential algebra
$\cQ^*_\infty$  to that $H^*(Y)$ of the fiber bundle $Y$.
\end{prop}

\begin{proof} The proof is obvious. The complex (\ref{lmp71}) is exact due to
the Poincar\'e lemma, and is a resolution of the constant sheaf $\Bbb R$ on
$J^\infty Y$ since $\gQ^*_\infty$ are sheaves of $\cQ^0_\infty$-modules.
Then, Theorem \ref{+132}  and  Lemma \ref{20jpa}
complete the proof.
\end{proof}

It follows that every closed form $\f\in \cQ^*_\infty$
splits into the sum
\beq
\f=\varphi +d\xi, \qquad \xi\in \cQ^*_\infty, \label{tams2} 
\eeq
where $\varphi$ is a closed form on the fiber bundle $Y$. This splitting 
plays an important role  in the sequel.
Since the graded differential algebras $\cO^*_\infty$ and $\cQ^*_\infty$
have the same de Rham cohomology, we further agree
to call 
\beq
H^*(J^\infty Y)\op=^{\rm def}H^*(\cQ^*_\infty)=H^*(\cO^*_\infty) \label{t43}
\eeq
the de Rham cohomology 
of $J^\infty Y$.

Let us consider the vertical exact sequence of sheaves
\beq
0\to \gO^m_X \ar^{\pi^{\infty*}} \gQ^{0,m}_\infty \ar^{d_V}\cdots \ar^{d_V} 
\gQ^{k,m}_\infty \ar^{d_V}\cdots, \qquad 0\leq m\leq n, \label{lmp90'}
\eeq
in the variational bicomplex
(\ref{7})
and the complex of their structure algebras
\beq
0\to \cO^m(X) \ar^{\pi^{\infty*}} \cQ^{0,m}_\infty \ar^{d_V}\cdots \ar^{d_V} 
\cQ^{k,m}_\infty \ar^{d_V}\cdots. \label{lmp90}
\eeq

\begin{prop} \label{7jp} 
There is an isomorphism 
\beq
H^*(m,d_V)=H^*(Y,\pi^*\gO^m_X) \label{lmp01}
\eeq
of the cohomology
groups $H^*(m,d_V)$ of the complex (\ref{lmp90}) to the cohomology groups
 $H^*(Y,\pi^*\gO^m_X)$ of $Y$ with
coefficients in the pull-back sheaf
$\pi^*\gO^m_X$ on $Y$.
\end{prop}

\begin{proof}
The exact sequence (\ref{lmp90'}) is a resolution of the pull-back sheaf
$\pi^{\infty*}\gO^m_X$ on $J^\infty Y$. Then, by virtue of Theorem \ref{+132},
we have a cohomology isomorphism
\be
H^*(m,d_V)=H^*(J^\infty Y,\pi^{\infty*}\gO^m_X). 
\ee
Lemma \ref{20jpa'} completes the proof.
\end{proof}

\begin{cor} Cohomology groups $H^{>\di Y}(m,d_V)$ vanish.
\end{cor}

The cohomology groups $H^*(m,d_V)$ have a
$C^\infty(X)$-module structure.  For instance, 
let 
\be
Y\cong X\times V\to X
\ee
be a trivial fiber bundle with a
typical fiber $V$. There is an obvious
isomorphism of $\Bbb R$-modules
\be
H^*(m,d_V)=\gO^m_X\ot H^*(V). 
\ee

Turn now to the rows of the variational bicomplex (\ref{7}). 
We have the exact sequence of sheaves 
\be
0\to \gQ^{k,0}_\infty \ar^{d_H}\gQ^{k,1}_\infty\ar^{d_H}\cdots  
\op\longrightarrow^{d_H} 
\gQ^{k,n}_\infty\ar^\tau\gE_k\to 0, \qquad k>0. 
\ee
Since the sheaves $\gQ^{k,0}_\infty$ and $\gE_k$ are fine, this is a
resolution of the fine sheaf $\gQ^{k,0}_\infty$. It states immediately the
following.

\begin{prop} \label{lmp99'}  The cohomology
groups $H^*(k,d_H)$ of the complex   
\beq
 0\to \cQ^{k,0}_\infty \ar^{d_H}\cQ^{k,1}_\infty\ar^{d_H}\cdots  
\op\longrightarrow^{d_H} 
\cQ^{k,n}_\infty\ar^{\tau} E_k\to 0, \qquad k>0, \label{am31}
\eeq
are trivial.
\end{prop}

This result at terms $\cQ^{k,<n}_\infty$ recovers that of \cite{tak2}. The
exactness of the complex (\ref{am31}) at the term $\cQ^{k,n}_\infty$ means
that, if
\be
\tau(\f)=0,
\qquad \f\in
\cQ^{k,n}_\infty,
\ee
then 
\be
\f=d_H\xi, \qquad \xi\in \cQ^{k,n-1}_\infty.
\ee
Since
$\tau$ is a projection operator, there is the $\Bbb R$-module decomposition
\beq
\cQ^{k,n}_\infty=E_k\oplus d_H(\cQ^{k,n-1}_\infty). \label{30jpa}
\eeq

\begin{rem}
One can derive Proposition \ref{lmp99'} from Theorem \ref{+132},
without appealing to that sheaves $\gE_k$ are acyclic. 
\end{rem}

Let us consider the exact sequence of sheaves
\be
0\to \Bbb R \to \gQ^0_\infty \ar^{d_H} \gQ^{0,1}_\infty
\ar^{d_H}\cdots  
\op\longrightarrow^{d_H} 
\gQ^{0,n}_\infty
\ee
where all sheaves except $\Bbb R$ are fine. Then, from Theorem \ref{+132} and
Lemma
\ref{20jpa}, we state the following.

\begin{prop} \label{lmp99} 
Cohomology groups $H^r(d_H)$, $r<n$, of the complex
\beq
0\to\Bbb R\to \cQ^0_\infty \ar^{d_H}\cQ^{0,1}_\infty\ar^{d_H}\cdots  
\op\longrightarrow^{d_H} \label{t10}
\cQ^{0,n}_\infty  
\eeq
are isomorphic to de Rham
cohomology groups $H^r(Y)$ of $Y$.
\end{prop}

This result recovers that of \cite{tak2}, but let us say something more.

\begin{prop} \label{lmp130} 
Any $d_H$-closed form $\si\in \cQ^{*,<n}_\infty$ is represented by the sum
\beq
\si= h_0\varphi + d_H\xi, \qquad \xi\in \cQ^*_\infty, \label{t40}
\eeq
where $\varphi\in \cO^{<n}_0$ is a closed form on the fiber bundle
$Y$. 
\end{prop}

\begin{proof}
Due to the relation 
\beq
 h_0d=d_Hh_0, \label{hn1}
\eeq
the horizontal projection $h_0$
provides a homomorphism of the de Rham complex (\ref{5.13'}) to the
complex 
\beq
0\to\Bbb R\to \cQ^0_\infty \ar^{d_H}\cQ^{0,1}_\infty\ar^{d_H}\cdots  
\op\longrightarrow^{d_H} 
\cQ^{0,n}_\infty\ar^{d_H} 0.   \label{+481'}
\eeq
Accordingly, there is a 
homomorphism
\beq
h_0^*: H^r(J^\infty Y) \to H^r(d_H), \qquad 0\leq r\leq n, \label{lmp125}
\eeq
of cohomology groups of these complexes. Proposition \ref{38jp} and
Proposition \ref{lmp99} show that, for $r< n$, the homomorphism (\ref{lmp125})
is an isomorphism (see the relation (\ref{1j}) below for the case
$r=n$). It follows that a horizontal form $\psi\in \cQ^{0,<n}$ is $d_H$-closed
(resp.
$d_H$-exact) if and only if $\psi=h_0\f$ where $\f$ is a closed (resp. exact)
form. The decomposition (\ref{tams2}) and Proposition
\ref{lmp99'} complete the proof.
\end{proof}

\begin{prop}
If
$\f\in\cQ^{0,<n}$ is a $d_H$-closed form, then $d_V\f=d\f$ is
necessarily $d_H$-exact.
\end{prop}

\begin{proof} 
Being nilpotent, the vertical differential $d_V$
defines a homomorphism of the complex (\ref{+481'}) to the complex
\be
 0\to \cQ^{1,0}_\infty \ar^{d_H}\cQ^{1,1}_\infty\ar^{d_H}\cdots  
\op\longrightarrow^{d_H} 
\cQ^{1,n}_\infty \ar^{d_H} 0
\ee
and, accordingly, a homomorphism of cohomology groups $H^*(d_H)\to
H^*(1,d_H)$ of these complexes. Since $H^{<n}(1,d_H)=0$, the result follows.
\end{proof}

Let us prolong the complex (\ref{t10}) to the variational complex
\beq
0\to\Bbb R\to \cQ^0_\infty \ar^{d_H}\cQ^{0,1}_\infty\ar^{d_H}\cdots  
\op\longrightarrow^{d_H} 
\cQ^{0,n}_\infty  \op\longrightarrow^\dl E_1 
\op\longrightarrow^\dl 
E_2 \longrightarrow \cdots\,  \label{b317}
\eeq
of the graded differential algebra $\cQ^*_\infty$.
In accordance with Lemma \ref{lmp03}, the variational complex (\ref{tams1})
is a resolution of the constant sheaf $\Bbb R$ on $J^\infty Y$.  Then,
Theorem \ref{+132} and Proposition
\ref{38jp} give immediately the following.

\begin{prop} \label{lmp05} 
There is an isomorphism
\beq
H^*_{\rm var}=H^*(Y) \label{lmp06}
\eeq
between cohomology $H^*_{\rm var}$ of the variational complex (\ref{b317}) and
de Rham cohomology of the fiber bundle $Y$.
\end{prop}

The isomorphism (\ref{lmp06}) recovers the result of \cite{tak2} and that
of \cite{ander80} at terms $\cQ^{0,n}_\infty$, $E_1$, but let us say
something more.
The relation (\ref{am13}) for $\tau$ and
the relation (\ref{hn1}) for $h_0$ define  a homomorphisms of the
de Rham complex (\ref{5.13'}) of the algebra $\cQ^*_\infty$ to the variational
complex (\ref{b317}). The corresponding homomorphism of their cohomology
groups is an isomorphism. Then, in accordance with the splitting
(\ref{tams2}), we come to the following assertion which complete Proposition
\ref{lmp130}.

\begin{prop} \label{t41} 
Any $\dl$-closed form $\si\in\cQ^{k,n}_\infty$, $k\geq 0$, is represented by
the sum
\ben
&& \si=h_0\varphi + d_H h_0\xi, \qquad k=0, \label{t42}\\
&& \si=\tau(\varphi) +\dl(\xi), \qquad k>0, \label{t42'}
\een
where $\varphi$ is a closed $(n+k)$-form on $Y$ and $\xi\in \cQ^*_\infty$. 
\end{prop}

\subs{5. Cohomology of $\cO^*_\infty$}

Thus, we have the whole cohomology of the graded
differential algebra $\cQ^*_\infty$. The following theorem provide us with
$d_H$- and $\dl$-cohomology of the graded
differential algebra
$\cO^*_\infty$. 

\begin{theo} \label{am11} 
Graded differential algebra $\cO^*_\infty$  has the same $d_H$- and
$\dl$-cohomology as $\cQ^*_\infty$.
\end{theo}

\begin{proof}
Let the common symbol $D$ stand for the coboundary operators  $d_H$ and
$\dl$ of the variational bicomplex. 
Bearing in mind the decompositions (\ref{t40}), (\ref{t42}) and (\ref{t42'}),
it suffices to show that, if an element
$\f\in
\cO^*_\infty$ is
$D$-exact with respect to the algebra $\cQ^*_\infty$ (i.e., $\f=D\varphi$,
$\varphi\in\cQ^*_\infty$), then it is
$D$-exact in the algebra
$\cO^*_\infty$ (i.e., $\f=D\varphi'$, $\varphi'\in\cO^*_\infty$).
Lemma \ref{am12} states that,
if
$Y$ is a contractible fiber bundle and a $D$-exact form $\f$ on $J^\infty Y$
is of finite jet order
$[\f]$ (i.e., $\f\in\cO^*_\infty$), there exists an exterior form $\varphi\in
\cO^*_\infty$ on $J^\infty Y$ such that $\f=D\varphi$. Moreover, a glance at
the homotopy operators for $d_V$, $d_H$ and $\dl$ \cite{olver} shows that  the
jet order
$[\varphi]$ of $\varphi$ is bounded for all exterior forms $\f$ of fixed
jet order. Let us call this fact the finite exactness of the operator
$D$. Given an arbitrary fiber bundle
$Y$, the finite exactness takes place on $J^\infty Y|_U$ over any open subset
$U$ of
$Y$ which is homeomorphic to a convex open subset of $\Bbb R^{\di Y}$.
Now, we show the following.

(i) Suppose that the finite exactness of the operator $D$ takes place on
$J^\infty Y$ over open subsets
$U$, $V$ of $Y$ and their non-empty overlap $U\cap V$. Then, it is also true on
$J^\infty Y|_{U\cup V}$.

(ii) Given a family $\{U_\al\}$ of disjoint open subsets of $Y$, let us
suppose that the finite exactness takes place on $J^\infty Y|_{U_\al}$ over
every subset $U_\al$ from this family. Then, it is true on $J^\infty Y$ over
the union
$\op\cup_\al U_\al$ of these subsets.

\noindent
If the assertions (i) and (ii) hold, the finite
exactness of
$D$ on $J^\infty Y$ takes place  since one can
construct the corresponding covering of the manifold $Y$
(\cite{bred2}, Lemma 9.5). 

{\it Proof of (i)}. Let
$\f=D\varphi\in\cO^*_\infty$ be a $D$-exact form on
$J^\infty Y$. By assumption, it can be brought into the form
$D\varphi_U$ on $(\pi^\infty_0)^{-1}(U)$ and $D\varphi_V$ on
$(\pi^\infty_0)^{-1}(V)$, where
$\varphi_U$ and $\varphi_V$ are exterior forms of finite jet
order. Due to the decompositions (\ref{t40}), (\ref{t42}) and (\ref{t42'}),
one can choose the forms $\varphi_U$, $\varphi_V$ such that
$\varphi-\varphi_U$ on $(\pi^\infty_0)^{-1}(U)$ and 
$\varphi-\varphi_V$ on $(\pi^\infty_0)^{-1}(V)$ are $D$-exact forms.
Let us consider the difference $\varphi_U-\varphi_V$ on 
$(\pi^\infty_0)^{-1}(U\cap V)$. It is a $D$-exact form of finite jet
order which, by assumption, can be written as 
$\varphi_U-\varphi_V=D\si$ where an exterior form
$\si$ is also of finite jet order. 
Lemma
\ref{am20} below shows that $\si=\si_U +\si_V$ where
$\si_U$ and
$\si_V$ are exterior forms of finite jet order on $(\pi^\infty_0)^{-1}(U)$ and
$(\pi^\infty_0)^{-1}(V)$, respectively. Then, putting
\be
\varphi'_U=\varphi_U-D\si_U, \qquad  
\varphi'_V=\varphi_V+ D\si_V,
\ee
we have the form $\f$  equal to
$D\varphi'_U$ on $(\pi^\infty_0)^{-1}(U)$ and
$D\varphi'_V$ on $(\pi^\infty_0)^{-1}(V)$, respectively. Since the
difference $\varphi'_U -\varphi'_V$ on $(\pi^\infty_0)^{-1}(U\cap V)$ vanishes,
we obtain $\f=D\varphi'$ on $(\pi^\infty_0)^{-1}(U\cup V)$ where 
\be
\varphi'\op=^{\rm def}\left\{
\begin{array}{ll}
\varphi'|_U=\varphi'_U, &\\
\varphi'|_V=\varphi'_V &
\end{array}\right.
\ee
is of finite jet order. 

{\it Proof of (ii)}. Let
$\f\in\cO^*_\infty$ be a $D$-exact form on
$J^\infty Y$.
The finite exactness on $(\pi^\infty_0)^{-1}(\cup
U_\al)$ holds since $\f=D\varphi_\al$ on every $(\pi^\infty_0)^{-1}(U_\al)$
and, as was mentioned above, the jet order
$[\varphi_\al]$ is bounded on the set of exterior forms
$D\varphi_\al$ of fixed jet order $[\f]$. 
\end{proof}

\begin{lem} \label{am20} 
Let $U$ and $V$ be open subsets of a fiber bundle $Y$ and $\si\in
\gO^*_\infty$ an exterior form of finite jet order on the non-empty overlap
$(\pi^\infty_0)^{-1}(U\cap V)\subset J^\infty Y$. Then, $\si$ splits
into  a sum $\si_U+ \si_V$ of exterior forms $\si_U$ and $\si_V$ of finite jet
order on
$(\pi^\infty_0)^{-1}(U)$ and $(\pi^\infty_0)^{-1}(V)$, respectively. 
\end{lem} 

\begin{proof}
By taking a smooth partition of unity on $U\cup V$ subordinate to the cover
$\{U,V\}$ and passing to the function with support in $V$, one gets a
smooth real function
$f$ on
$U\cup V$ which is 0 on a neighborhood of $U-V$ and 1 on a neighborhood of
$V-U$ in $U\cup V$. Let $(\pi^\infty_0)^*f$ be the pull-back of $f$ onto
$(\pi^\infty_0)^{-1}(U\cup V)$. The exterior form $((\pi^\infty_0)^*f)\si$ is
zero on a neighborhood of $(\pi^\infty_0)^{-1}(U)$ and, therefore, can be
extended by 0 to $(\pi^\infty_0)^{-1}(U)$. Let us denote it $\si_U$.
Accordingly, the exterior form
$(1-(\pi^\infty_0)^*f)\si$ has an extension $\si_V$ by 0 to 
$(\pi^\infty_0)^{-1}(V)$. Then, $\si=\si_U +\si_V$ is a desired decomposition
because $\si_U$ and $\si_V$
are of finite jet order which does not exceed that of $\si$. 
\end{proof}

It is readily observed that Theorem \ref{am11} can be applied to de
Rham cohomology of $\cO^*_\infty$ whose isomorphism (\ref{t43}) to that of
$\cQ^*_\infty$ has been stated. 

\subs{6. The global inverse problem}

The variational complex (\ref{b317}) provides the algebraic approach to the
calculus of variations on fiber bundles in the class of
exterior forms of locally finite jet order \cite{bau,book,tul}. For instance,
the variational operator 
$\dl$ acting on
$\cQ^{0,n}_\infty$ is the Euler--Lagrange map, while $\dl$ acting on $E_1$ is
the Helmholtz--Sonin map.  Let $L\in \cQ^{0,n}_\infty$
be a horizontal density on $J^\infty Y$. One can think of $L$ as
being a Lagrangian of local finite order.
Then, the 
decomposition (\ref{30jpa}) leads to the first variational formula
\be
dL=\tau(dL) + (\id -\tau)(dL)= \dl(L) + d_H(\phi), \qquad \phi\in
\cQ^{1,n-1}_\infty, 
\ee
where $\dl (L)$ 
is the Euler--Lagrange form associated with the Lagrangian $L$.

Let us relate the cohomology isomorphism (\ref{lmp06}) to the global inverse
problem of the calculus of variations. As a particular repetition of
Proposition \ref{t41}, we come to its following solution in the class of
exterior forms of locally finite jet order.

\begin{theo} \label{lmp112'} 
A Lagrangian $L\in \cQ^{0,n}_\infty$ is variationally trivial, i.e., 
$\dl(L)=0$ if and only if 
\beq
L=h_0\varphi + d_H h_0\xi, \qquad \xi\in \cQ^*_\infty, \label{tams3}
\eeq
where $\varphi$ is a closed $n$-form on $Y$ (see the expression (\ref{t42})). 
\end{theo}

\begin{theo} \label{lmp113'} 
An 
Euler--Lagrange-type operator $\cE\in E_1$ satisfies the Helmholtz
condition $\dl(\cE)=0$ if and only if 
\be
\cE=\dl(L) + \tau(\f), \qquad L\in\cQ^{0,n}_\infty, 
\ee
where $\f$ is a closed $(n+1)$-form on $Y$ (see the expression (\ref{t42'})).
\end{theo}

Theorem \ref{lmp113'} contains the similar result of \cite{ander80,tak2}.

\begin{rem}
As a consequence of Theorem \ref{lmp112'}, one obtains that 
the cohomology group
$H^n(d_H)$ of the complex (\ref{+481'}) obeys the relation
\beq
H^n(d_H)/H^n(Y)=\dl(\cQ^{0,n}_\infty), \label{1j}
\eeq 
where $\dl(\cQ^{0,n}_\infty)$ is the $\Bbb R$-module of Euler--Lagrange
forms on $J^\infty Y$.
\end{rem}

Theorem \ref{am11} provides the similar solution of the global inverse
problem in the class of exterior forms of finite jet order.
The theses of Theorem \ref{lmp112'} and Theorem \ref{lmp113'} remain true
if all exterior forms belong to $\cO^*_\infty$.
Theorem \ref{lmp113'} contains the result of \cite{vin}. 

As was mentioned above, a solution of the global inverse problem in the
calculus of variations in the class of exterior forms of a fixed jet order has
been suggested in
\cite{ander80} by a computation of cohomology of the fixed
order variational sequence.
The first thesis of \cite{ander80} agrees with Theorem \ref{lmp112'} for
finite order Lagrangians, but says that the
jet order of the form $\xi$ in the expression (\ref{tams3}) is $k-1$ if $L$ is
a $k$-order variationally trivial Lagrangian. The second one states that a
$2k$-order Euler--Lagrange operator can be always  associated with a $k$-order
Lagrangian. However, because of the sophisticated technique, these results
were not widely recognized.

\bigskip\bigskip

{\parindent=0pt

{\sc Department of Mathematics and Physics, University of Camerino, 62032
Camerino (MC), Italy}

{\it E-mail address}: {\bf giachetta@campus.unicam.it}
\bigskip

{\sc Department of Mathematics and Physics, University of Camerino, 62032
Camerino (MC), Italy}

{\it E-mail address}: {\bf mangiaro@camserv.unicam.it}
\bigskip

{\sc Department of Theoretical Physics, Moscow State University, 117234
Moscow, Russia}

{\it E-mail address}: {\bf sard@grav.phys.msu.su}

}

\bigskip
\bigskip

\begin{thebibliography}{aaa}

\bibitem{ander80} {\sc I. Anderson and T. Duchamp}, On the existence of global
variational principles, Amer. J. Math., {\bf 102} (1980), 781-868. MR {\bf
82d}:58027

\bibitem{ander} {\sc I. Anderson}, Introduction to the variational bicomplex,
Contemp. Math., {\bf 132} (1992), 51-74. MR {\bf 94b}: 58045

\bibitem{bau} {\sc M. Bauderon}, Differential geometry and Lagrangian formalism
in the calculus of variations, in {\sl Differential Geometry, Calculus of
Variations, and their Applications}, (G. Rassias and T. Rassias, eds), 67-82,
Lecture Notes in Pure and Applied Mathematics {\bf 100}, Marcel Dekker Inc.,
New York, 1985. MR {\bf 87h}: 58050

\bibitem{bred} {\sc G. Bredon}, {\sl Sheaf Theory},
McGraw-Hill Book Company, New York, 1967. MR {\bf 36}:4570

\bibitem{bred2} \underline{\hskip1.5cm}, {\sl Topology and Geometry},
Graduate Texts in Mathematics {\bf 139}, Springer-Verlag, Berlin, 1997. MR
{\bf 2000b}:55001


\bibitem{hir} {\sc F. Hirzebruch}, {\sl Topological Methods in Algebraic
Geometry}, Springer-Verlag, Berlin, 1966. MR {\bf 34}:2573

\bibitem{book} {\sc G. Giachetta, L. Mangiarotti and G. Sardanashvily}, {\sl
New Lagrangian and Hamiltonian Methods in Field Theory}, World Scientific,
Singapore, 1997.

\bibitem{god} {\sc R. Godement}, {\sl Th\'eorie des faisceaux}, Hermann,
Paris, 1964. 

\bibitem{kash} {\sc M. Kashiwara and P. Scapira}, {\sl Sheaves on Manifolds},
A Series of Comprehensive Studies in Mathematics, {\bf 292}, Springer-Verlag,
Berlin, 1990. MR {\bf 92a}:58132

\bibitem{kru90} {\sc D. Krupka}, Variational sequences on finite order jet
spaces, in {\sl Proceeding of the Conference on Differential Geometry and its
Applications (Brno, 1989)}, 236-254, World Scientific, Singapore, 1990. MR {\bf
91i}:58037

\bibitem{kru98} {\sc D. Krupka and J. Musilova}, Trivial Lagrangians in field
theory, Diff. Geom. Appl., {\bf 9} (1998), 293-305. MR {\bf 2000g}:58027a

\bibitem{man} {\sc L. Mangiarotti and M. Modugno}, Graded Lie algebras and
connections on a fibered space, J. Math. Pures et Appl., {\bf 63} (1984),
111-120. MR {\bf 86c}:53014

\bibitem{olver} {\sc P. Olver}, {\sl Applications of Lie Groups to
Differential Equations}, Springer-Verlag, Berlin, 1997. MR {\bf 94g}:58260

\bibitem{sard} {\sc G. Sardanashvily}, {\sl Gauge Theory in Jet Manifolds},
Hadronic Press Monographs in Applied Mathematics, Hadronic Press Inc., Palm
Harbor, Fl, 1993. MR {\bf 95e}:58054

\bibitem{tak1} {\sc F. Takens}, Symmetries, conservation laws and variational
principles, in {\sl Geometry and Topology}, (J. Palis and M. do Carmo, eds),
581-604, Lect. Notes in Mathematics {\bf 597}, Springer-Verlag, Berlin, 1977.
MR {\bf 50}:31255

\bibitem{tak2} \underline{\hskip1.5cm}, A global version of the inverse problem
of the calculus of variations, J. Diff. Geom., {\bf 14} (1979), 543-562. MR
{\bf 83b}:58028

\bibitem{tul} {\sc W. Tulczyjew}, The Euler--Lagrange resolution, 
in {\sl Differential
Geometric Methods in Mathematical Physics},(P. Garc\'\i a, A. P\'erez-Rend\'on
and J.Souriau, eds) 22-48, Lect. Notes in Mathematics {\bf 836},
Springer-Verlag, Berlin, 1980. MR {\bf 82j}:58034

\bibitem{vin} {\sc A. Vinogradov} The ${\cal C}$-spectral sequence, Lagrangian
formalism, and conservation laws. II. The nonlinear theory., J. Math. Anal.
Appl., {\bf 100} (1984), 41-129. MR {\bf 85j}:58150b

\bibitem{vit} {\sc R. Vitolo}, Finite order variational bicomplex, Math.
Proc. Cambridge Phil. Soc., {\bf 125} (1998), 321-333. MR {\bf 99j}:58232

\end{thebibliography}
\end{document}